\documentclass[11pt]{article}
\usepackage[utf8]{inputenc}
\usepackage[T1]{fontenc}

\usepackage{amsfonts}
\usepackage{amsmath}
\usepackage{amssymb}
\usepackage{amsthm}

\usepackage{epigraph}

\usepackage{xcolor}
\usepackage{hyperref}

\usepackage[shortlabels]{enumitem}

\theoremstyle{definition}
\newtheorem{tw}{Theorem}

\newtheorem*{xp}{Example}
\newtheorem{df}{Definition}
\newtheorem{pr}{Proposition}

\newtheorem*{con}{Conjecture}

\DeclareMathOperator\gen{gen}

\newcommand{\vphi}{\varphi}
\newcommand{\mtt}{\mathsf}
\newcommand{\mfr}{\mathfrak}

\newcommand{\msf}{\mathsf}

\newcommand{\eps}{\varepsilon}
\newcommand{\lf}{\left \lfloor}
\newcommand{\rf}{\right \rfloor}

\providecommand{\keywords}[1]
{
  \small	
  \textbf{\textit{Keywords:}} #1
}

\title{Finding binary words with a given number of subsequences}
\author{Radosław Żak\footnote{Undergraduate student;  Jagiellonian University, 
Faculty of Mathematics and Computer Science; 
Łojasiewicza 6, 30-348 Kraków, Poland.\newline \texttt{radoslaw.zak@student.uj.edu.pl}}}
\date{February 18, 2022}

\begin{document}

\maketitle

\begin{abstract}
We relate binary words with a given number of subsequences to continued fractions of rational numbers with a given denominator. We deduce that there are binary strings of length $O(\log n \log \log n)$  with exactly $n$ subsequences; this can be improved to $O(\log n)$ under assumption of Zaremba's conjecture.
\end{abstract}


\section{Introduction}

The number of subsequences of a binary word was investigated so far mainly from the probabilistic point of view. Collins in \cite{collins} proved that a random binary string of length $n$ has $2\left(\frac{3}{2}\right)^n-1$ subsequences on average, while Biers-Ariel, Godbole and Kelley in \cite{ariel} generalized this result to the case where the probabilities of occurrence of particular letters are distinct. Flaxman, Harrow and Sorkin proved in \cite{flaxman} that the string maximizing the number of subsequences consists of cyclically repeating letters.

Here our main subject of interest is the problem of finding a binary word with exactly $n$ subsequences, as short as possible. We exhibit a relation between binary strings, Euclidean algorithm and continued fractions. Using results of Rukavishnikova (\cite{rukavishnikova}), we can deduce that there are words with $n$ subsequences and length $O(\log n \log \log n)$. If Zaremba's conjecture is true, we can even find such a word with length $\Theta(\log n)$.

We also derive two interesting facts. Theorem \ref{add} states that if we have a word $\mfr{s}$ and we want to add $n$ letters at the end, so that the resulting string has as many subsequences as possible, then $\msf{ABAB}\ldots$ or $\msf{BABA}\ldots$ is an optimal choice for the letters added (depending on the last letter of $\mfr{s}$). Theorem \ref{approx} relates good approximations of an irrational number with restricted partial quotients by a rational with denominator $N$ to short words with $N-1$ subsequences.

For a clear distinction between letters, numbers and words we use alphabet $\{ \msf{A},\, \msf{B}\}$, and we denote words with Fraktur, eg. as $\mfr{s},\, \mfr{t}$. We also use the following notations:

\begin{itemize}
    \item $\msf{0}$ -- the empty string.
    \item $\mfr{s} \circ \mfr{t}$ -- the concatenation of $\mfr{s}$ and $\mfr{t}$.
    \item $\mfr{s}^k$ is simply $\mfr{s} \circ \mfr{s} \circ \ldots \circ \mfr{s}$ -- $k$ copies of $\mfr{s}$ concatenated.
    \item $\msf{P}(\mfr{s})$ -- the number of subsequences of $\mfr{s}$.
    \item $\mfr{s}^*$ -- the word created by replacing in $\mfr{s}$ all letters $\msf{A}$ with $\msf{B}$ and \textit{vice versa}.
    \item $|\mfr{s}|$ -- the length of string $\mfr{s}$.
    \item $\vphi(\cdot)$ -- the Euler function.
\end{itemize}

\section{Words and the Euclidean algorithm}

The following notion characterizes words with a given number of subsequences in a surprising way.

\begin{df}
For two coprime integers $a,\, b \geqslant 1$ we construct the word $\gen(a,b)$ recursively:

\[ \gen(a,b) = \begin{cases} \msf{0} \textrm{ if } a=b=1,\\
\msf{A}\circ \gen(a-b,b)\ \textrm{ if } a>b,\\
\msf{B}\circ \gen(a,b-a)\ \textrm{ if } b>a.
\end{cases} \]
\end{df}

This is, in some sense, the description of the Euclidean algorithm for numbers $a,\, b$ (in particular this algorithm implies that words $\gen(a,b)$ are well-defined). As it turns out, we have a large control over number of subsequences of this word.

\begin{tw} \label{twGen}
$\msf{P}(\gen(a,b))=a+b-1.$
\end{tw}

Before proving the theorem, let us introduce an additional notation.

\begin{df}
If $\mfr{s}$ is a binary word, denote as $\msf{P}^{\msf{A}}(\mfr{s})$ the number of subsequences of $\mfr{s}$ starting with $\msf{A}$, \textbf{including the empty one}, and as $\msf{P}_\msf{A}$ the number of subsequences ending with $\msf{A}$ (empty string also included). Similarly, $\msf{P}^{\msf{B}}(\mfr{s})$ and $\msf{P}_{\msf{B}}(\mfr{s})$ are the numbers of subsequences of $\mfr{s}$ respectively starting and ending with $\msf{B}$.
\end{df}

This way $\msf{P}(\mfr{s}) = \msf{P}^\msf{A}(\mfr{s})+\msf{P}^\msf{B}(\mfr{s})-1$ (the $-1$ corresponds to the empty subsequence counted two times). 

\begin{proof}
We will be proving inductively that $\msf{P}^{\msf{A}}(\gen(a,b))=a$. Then the claim follows easily, as similarly we get $\msf{P}^{\msf{B}}(\gen(a,b))=b$. For $a=b=1$ the empty word is the only subsequence of $\gen(1,1)$.

If $a<b$, the first letter of $\gen(a,b)$ is $\msf{B}$. It cannot be contained in any subsequence starting with $\msf{A}$, so deleting it does not change the value of $\msf{P^A}$. But $\gen(a,b) = \msf{B} \circ \gen(a,b-a)$, and by induction assumption $\gen(a,b-a)$ has exactly $a$ subsequences starting with $\msf{A}$.

Now consider case $a>b$. The first letter of $\gen(a,b)$ is $\msf{A}$. Observe that if we choose letters from $\gen(a,b)$ in order to form a subsequence starting with $\msf{A}$, we may use the first letter as a start -- the only exception is the empty sequence. Therefore, as $\gen(a,b) = \msf{A} \circ \gen(a-b,b)$ we have a correspondence between nonempty subsequences of $\gen(a,b)$ starting with $\msf{A}$ and subsequences of $\gen(a-b,b)$. By induction assumption there are $a-1$ of the latter, which along with the empty word form $a$ subsequences of $\gen(a,b)$ starting with $\msf{A}$. This ends the proof.

\end{proof}

\begin{xp}
Take $a=11,\, b=7$. Then $\gen(11,7)=\msf{A} \circ \gen(4,7) = \msf{AB} \circ \gen(4,3) = \msf{ABA} \circ \gen(1,3) = \msf{ABAB} \circ \gen(1,2) = \msf{ABABB}$. We can check that $\msf{ABABB}$ has $10$ nonempty subsequences starting with $\msf{A}$, and $6$ starting with $\msf{B}$.
\end{xp}

The proof of Theorem \ref{twGen} says that the inverse function of $\gen$ is $\mfr{s} \to \left( \msf{P^A}(\mfr{s}), \msf{P^B}(\mfr{s}) \right)$. As it is not hard to see that any binary word can be written as $\gen(a,b)$ for some $a,\, b$, this gives a bijection between binary words and pairs of coprime positive integers, which allows us to form the following corollary (insignificant for our later reasonings, but interesting on its own).

\begin{pr}
There are exactly $\vphi(N+1)$ binary words with $N$ subsequences.
\end{pr}

\begin{proof}
All of them have form $\gen(a,b)$, where $b=N+1-a$ and $\gcd(a,b)=1$. The latter is equivalent to $\gcd(a,N+1)=1$, so there are exactly $\vphi(N+1)$ possible choices for $a$.
\end{proof}

\subsection{Concatenation theorem}

Before we proceed further, it is important to state the formula for the number of subsequences of two words' concatenation.

\begin{tw} \label{twKonKat}
\[ \mtt{P}(\mfr{s} \circ \mfr{t}) = \msf{P_A}(\mfr{s}) \msf{P^B}(\mfr{t}) + \msf{P_B}(\mfr{s}) \msf{P^A}(\mfr{t}) - 1. \]
\end{tw}

\begin{proof}
We prove inductively on $|\mfr{s}|$. For $\mfr{s}=\msf{0}$ the claim is clear. Suppose $\mfr{s}$ is nonempty. Without loss of generality $\mfr{s} = \mfr{u} \circ \msf{A}$ for some $\mfr{u}$. Then $\msf{P_A}(\mfr{s}) = \msf{P_A}(\mfr{u}) + \msf{P_B}(\mfr{u})$ and $\msf{P_B}(\mfr{s}) = \msf{P_B}(\mfr{u})$. Therefore:

\begin{align*}
    \mtt{P}_{\msf{A}}(\mfr{s}) \mtt{P}^{\msf{B}}(\mfr{t}) + \mtt{P}_{\msf{B}}(\mfr{s}) \mtt{P}^{\msf{A}}(\mfr{t}) =\\
    \left(\mtt{P}_{\msf{A}}(\mfr{u})+\mtt{P}_{\msf{B}}(\mfr{u}) \right) \mtt{P}^{\msf{B}}(\mfr{t}) + \mtt{P}_{\msf{B}}(\mfr{u}) \mtt{P}^{\msf{A}}(\mfr{t}) = \\
    \mtt{P}_{\msf{A}}(\mfr{u}) \mtt{P}^{\msf{B}}(\mfr{t}) + \mtt{P}_{\msf{B}}(\mfr{u}) \left(\mtt{P}^{\msf{B}}(\mfr{t})+\mtt{P}^{\msf{A}}(\mfr{t}) \right) = \\
    \mtt{P}_{\msf{A}}(\mfr{u}) \mtt{P}^{\msf{B}}(\msf{A} \circ \mfr{t}) + \mtt{P}_{\msf{B}}(\mfr{u}) \mtt{P}^{\msf{A}}(\msf{A} \circ \mfr{t}).
\end{align*}
\end{proof}


\section{Short words}

Define string $\mfr{z}_n$ by $\mfr{z}_0 = \msf{0}$, $\mfr{z}_{n+1} = \msf{A} \circ \mfr{z}^*_n$. In other words, $\mfr{z}_n$ is a string $\msf{ABABA} \ldots$ with $n$ letters.

\begin{tw} \label{add}
Let $\mfr{s}$ be a word ending with $\msf{B}$ (or an empty one), and $n$ a non-negative integer. Then

\[ \max_{|\mfr{t}|=n} \msf{P}(\mfr{s} \circ \mfr{t}) = \msf{P}(\mfr{s} \circ \mfr{z}_n). \]

Moreover if $\mfr{s} \neq \msf{0}$, the maximum is attained only by $\mfr{z}_n$.
\end{tw}

\begin{proof}
We will be proving inductively on $n$. For $n=0$ it is clear. Suppose $n>0$ and let $\mfr{t}$ be the $n$-letter word for which $\msf{P}(\mfr{s} \circ \mfr{t})$ is maximal. When $\mfr{t}$ starts with $\msf{A}$, we get the claim by applying inductive hypothesis to $(\mfr{s} \circ \msf{A})^*$ and $n-1$. When $\mfr{t}$ starts with $\msf{B}$, by applying inductive hypothesis to $\mfr{s} \circ \msf{B}$ and $n-1$ we can replace $\mfr{t}$ with $\mfr{z}_n^*$. Hence we just need to check that $\msf{P}(\mfr{s} \circ \mfr{z}_n) \geqslant \msf{P}(\mfr{s} \circ \mfr{z}^*_n)$ (and that the inequality is strict whenever $\mfr{s}$ is nonempty).

If $\mfr{s}=\msf{0}$, this is clear. Suppose $|\mfr{s}|>0$. By Theorem \ref{twKonKat} we know that $\msf{P}(\mfr{s} \circ \mfr{t})+1 = \msf{P_A}(\mfr{s}) \msf{P^B}(\mfr{t}) + \msf{P_B} (\mfr{s}) \msf{P^A}(\mfr{t})$. We know that $\msf{P^A}(\mfr{z}_n)=\msf{P^B}(\mfr{z}^*_n)$ and $\msf{P^B}(\mfr{z}_n)=\msf{P^A}(\mfr{z}^*_n)$. Since $\mfr{s}$ ends with $\msf{B}$, we have $\msf{P_B}(\mfr{s}) > \msf{P_A} (\mfr{s})$. Moreover $\msf{P^A}(\mfr{z}_n) > \msf{P^A}(\mfr{z}^*_n)$, therefore by rearrangement inequality $\msf{P}(\mfr{s} \circ \mfr{z}_n) > \msf{P}(\mfr{s} \circ \mfr{z}^*_n)$.
\end{proof}

This theorem allows us to replicate the result from \cite{flaxman} for binary strings:

\begin{pr}
The words $\mfr{z}_n$ and $\mfr{z}_n^*$ (and only them) have the maximal number of subsequences among binary words on $n$ letters.
\end{pr}

\begin{proof}
For $n=0$ it is clear; for $n>0$ by using Theorem \ref{add} for $\mfr{s}=\msf{B}$ we get that $\mfr{z}_n^*$ has more subsequences than any string starting with $\msf{B}$; for $\mfr{z}_n$ it is symmetric.
\end{proof}

Using induction we can enumerate $\msf{P}(\mfr{z}_n) = F_{n+3}-1$, where $F_m$ is the $m$-th Fibonacci number. Indeed,
\[ \msf{P}(\mfr{z}_n)+1= \msf{P}^\msf{A}(\msf{A}\circ \mfr{z}^*_{n-1}) + \msf{P}^\msf{B}(\msf{AB} \circ \mfr{z}_{n-2}) = (1+ \msf{P}(\mfr{z}^*_{n-1})) + (1+ \msf{P}(\mfr{z}_{n-2})) \]

(being careful for the empty subsequence), which by inductive assumption is equal to $F_{n+2}+F_{n+1}=F_{n+3}$.

Since $F_n \approx \frac{1}{\sqrt{5}} \phi^n$, where $\phi=\frac{1+\sqrt{5}}{2}$ is the golden ratio, we get that a word with $n$ subsequences has to have length $\Omega(\log n)$. The following conjecture seems natural:

\begin{con}
For any $n \geqslant 1$ there is a binary word with length $O(\log n)$ and exactly $n$ subsequences.
\end{con}

\subsection{Continued fractions}

To arrive at our main results, we observe the following duality:

\begin{pr} \label{cont}
Let $a,\, b$ be coprime positive integers and $\theta=\frac{a}{b}$. If the continued fraction of $\theta$ is $[c_0;c_1,c_2,\ldots,c_k+1]$, then $\gen(a,b)=\msf{A}^{c_0}\msf{B}^{c_1}\msf{A}^{c_2}\ldots$.
\end{pr}

\begin{proof}
We proceed by induction. If $a=b=1$, then $\theta=1$ and $\gen(a,b) = \msf{A}^0 = \msf{0}$. If $b>a$, $c_0=0$. We can swap these numbers, so that $\theta$ becomes $\theta^{-1} = [c_1;c_2,c_3,\ldots]$ and proceed further. Finally if $a>b$, let $a':=a-b$. Then $\gen(a,b) = \msf{A} \circ \gen(a',b)$, and $\frac{a'}{b} = \frac{a}{b} - 1 = [c_0-1; c_1,c_2,\ldots,c_k+1]$. This finishes the proof.
\end{proof}

It is tempting to look for short strings with $N-1$ subsequences among those words which start with long segment of form $\mfr{z}_k$ -- after all, we know this can give us a lot of subsequences. For that, take $a \approx \frac{1}{\phi} N$ and $b=N-a$; then $\frac{a}{b} \approx \phi$. Theorem \ref{approx} captures the same idea for an arbitrary irrational number, such that the terms of its continued fraction are bounded.

Choose an integer $C \geqslant 1$. We say that a number $\theta$ has partial quotients bounded by $C$ if the continued fraction of $\theta=[c_0;c_1,c_2,\ldots]$ has $c_i\leqslant C$ for all $i$. Denote by $S_C$ the set of irrational numbers with partial quotients bounded by $C$. It is clearly closed under operation $\theta \mapsto \theta^{-1}$, and if $\theta \in S_C$, $\theta > 1$, then also $\theta-1 \in S_C$.

\begin{tw} \label{approx}
Take $\xi \in (0,1) \cap S_C$. If $\xi$ can be approximated by an irreducible fraction $\frac{a}{N}$ with error $\left \vert \xi - \frac{a}{N} \right \vert = \delta$, then $\gen(a,N-a)$ is a binary word with $N-1$ subsequences and length
\[ O(C \log N + N \sqrt{\delta C^3}). \]
\end{tw}

\begin{proof}
Let $\theta = \frac{1}{\xi} - 1$. Since $\xi<1$, $\theta$ is positive, moreover $\xi \in S_C$ implies $\theta \in S_C$.

Let $b=N-a$, $x= \xi N$, $y = (1-\xi) N$. Now $\frac{y}{x} = \theta$, $x+y=N$ and $\eps := |y-b|= |x-a| = N |\xi - \frac{a}{N}| = N\delta$.

We can apply Euclidean algorithm to pairs $(x,y)$ and $(a,b)$. They go the same way, until the total error $|x-a|+|y-b|$ is greater than (or equal to) $|x-y|$ (when $|x-a|+|y-b|<|x-y|$ the equivalence $a<b \iff x<y$ holds). Let $\mfr{s}$ be the part of $\gen(a,b)$ corresponding to that interval of time. 

Observe that the error $(x-a,y-b)$ starts at $(\pm \eps, \mp \eps)$, and while we execute the Euclidean algorithm, taking a difference between $x$ and $y$ corresponds to adding absolute value of one of these errors to another. Therefore the errors after we have the word $\mfr{s}$ are $(\pm \msf{P}_\msf{A}(\mfr{s}) \eps,\mp \msf{P}_\msf{B}(\mfr{s}) \eps)$, and the sum of their absolute values is $(\msf{P}(\mfr{s})+1) \eps$.

Let $(z,w)$ be the values of $(x,y)$ after we execute on them part of the Euclidean algorithm corresponding to $\mfr{s}$. Since $\left(\max(z,w),\max(z,w)\right) \geqslant (z,w)$ (coordinate-wise), by reversing the algorithm (adding one coordinate to another, the choice of coordinate is indicated by $\mfr{s}$), which preserves inequalities (it consists only of adding) we get $\left(\max(z,w) \msf{P}^\msf{A}(\mfr{s}),\max(z,w) \msf{P}^\msf{B}(\mfr{s})\right) \geqslant (x,y)$, thus $\max(z,w) \geqslant \frac{x+y}{\msf{P}(\mfr{s})+1} = \frac{N}{\msf{P}(\mfr{s})+1}$.

The value of $\frac{\max(x,y)}{\min(x,y)}$ in one step of the Euclidean algoritm changes either by $\theta \to \theta-1$ or $\theta \to \frac{1}{\theta-1}$. In both cases it stays inside $S_C$. Now

\[ |z-w| = \left(1- \frac{\min(z,w)}{\max(z,w)}\right) \max(z,w) \]

and the fraction lies in $(0,1) \cap S_C$ ($z \neq w$, since $\xi$ is irrational). The maximum of $(0,1) \cap S_C$ is not greater than $\frac{1}{1+\frac{1}{C+\frac{1}{1}}} = \frac{C+1}{C+2}$, so

\[ |z-w| \geqslant \frac{\max(z,w)}{C+2} \geqslant \frac{N}{(C+2)(\msf{P}(\mfr{s})+1)}. \]

Since $|z-w|$ is not greater than the total error $\eps (\msf{P}(s)+1)$ (if it was, we could continue the algorithm), we obtain

\[ \msf{P}(s)+1 \geqslant \sqrt{\frac{N}{(C+2)\eps}}. \]

We can write $\gen(a,b)$ as $\mfr{s} \circ \mfr{t}$ for some word $\mfr{t}$. Let us say that $\mfr{t}$ starts with $\msf{B}$ (the other case is analogous). By Theorem \ref{twKonKat},

\[ N = \msf{P}(\gen(a,b))+1 \geqslant \msf{P}_\msf{A}(\mfr{s}) \msf{P}^\msf{B}(\mfr{t}) \geqslant \msf{P}_\msf{A}(\mfr{s}) |\mfr{t}| \]

(the last inequality follows from the fact that we have subsequences of $\mfr{t}$ starting with $\msf{B}$ with all possible lengths).

Since $\mfr{s}$ comes from the continued fraction of $\theta$, which lies in $S_C$, each letter can repeat at most $C$ times in a row. By easy induction, if $\mfr{w}$ is any word with exactly $k$ letters $\msf{B}$ at the end, then $\msf{P}^\msf{A}(\mfr{w}) \geqslant \frac{1}{k+2} (\msf{P}(\mfr{w})+1)$ (equivalently $\msf{P}^\msf{B}(\mfr{w}) \leqslant (k+1)\msf{P}^\msf{A}(\mfr{w})$; for $k=0$ it works, and step $k \to k+1$ changes both sides of the equation by $\msf{P}^\msf{A}(\mfr{w})$). Therefore

\[ N \geqslant \msf{P}_\msf{A}(\mfr{s}) |\mfr{t}| \geqslant \frac{1}{C+2} (\msf{P}(\mfr{s})+1) |\mfr{t}| \geqslant \sqrt{\frac{N}{(C+2)^3 \eps}} |\mfr{t}|, \]

and so, putting $\eps = N \delta$,

\[ N \sqrt{\delta (C+2)^3} \geqslant |\mfr{t}|. \]

Now it is left to see that $|\mfr{s}| = O(C \log N)$. Indeed, every letter in $\mfr{s}$ comes at most $C$ times in a row, so $\mfr{s}$ has $\mfr{z}_m$ or $\mfr{z}^*_m$ as a subsequence, where $m=\lf \frac{1}{C} |\mfr{s}| \rf$. On the other hand, then $\msf{P}(\mfr{z}_m) \leqslant \msf{P}(\mfr{s}) \leqslant N$, so $m=O(\log N)$.

\end{proof}

The second summand, seemingly linear in $N$, is actually balanced by $\delta$. If $\delta = O(N^{-2})$ (which is the best possible option), it reduces to $\sqrt{C^3}$, which is even smaller than the $C \log N$ when $C$ does not exceed $\log^2 N$. Unluckily without some results on diophatine approximation we cannot give any particular bounds on how small $\delta$ can be.

However, the case of rational $\xi$ seems to be better studied in theory of continued fractions. By Proposition \ref{cont}, now the correspondence between fractions and words is even clearer. Zaremba conjectured in 1972 that for any $N$ there is an irreducible fraction $\frac{a}{N}$ that has all partial quotients bounded by $5$. In our case, this would imply that $\gen(a,N-a)$ is a word with $N-1$ subsequences, and length logarithmic in $N$. Bourgain and Kontorovich (\cite{bourgain}) proved that Zaremba's hypothesis is true (with bound $C=50$) for $N$ forming a set of density 1 in $\mathbb{N}$.

If we denote by $S_N(a)$ the sum of partial quotients in $\frac{a}{N}$, by Proposition \ref{cont} we have $S_N(a)=|\gen(a,N-a)|+2$. Rukavishnikova proved in \cite{rukavishnikova} an analog of the law of large numbers for $S_N(a)$, namely that if $g(N)$ is any unboundedly increasing function, such that $g(N) \leqslant \sqrt{\log \log N}$, then the fraction of all numbers $S_N(a)$ that fall outside of the interval $|S_N(a)-\frac{12}{\pi^2} \log N \log \log N| \leqslant g(N) \log N \sqrt{\log \log N}$ grows asympotically slower than $\frac{1}{g(N)^2}$. Thus we should expect that for most values of $N$, median of all numbers $S_N(a)$ is around $\frac{12}{\pi^2} \log N \log \log N$. Unfortunately, if $N$ has many prime divisors, then $\# \mathbb{Z}^*_N$ can be smaller than $\frac{N}{g(N)^2}$, so this result mean nothing about values of $S_N(a)$ with $\gcd(a,N)=1$ for general $N$ -- the case we are interested in.

Thus we need to use another result by Rukavishnikova, also featured in $\cite{rukavishnikova}$:

\begin{tw}[Rukavishnikova]
Suppose that $g(d)$ is unboundedly increasing sequence of positive real numbers for which $g(d) \leqslant (\log d)^2$. Then for $d>2$

\[ \frac{1}{\vphi(d)} \# \{a \in \mathbb{Z}_d^* : S_d(a) \geqslant g(d) \log d \log \log d \} = O\left( \frac{1}{g(d)}\right). \]
\end{tw}

In particular, the smallest value of $S_N(a)$ among $a \in \mathbb{Z}^*_N$ is $O(\log N \log \log N)$ (as it is growing slower than any function growing faster than $\log N \log \log N$). This implies the following:

\begin{pr}
For any positive integer $N$ there is a binary word of length $O(\log N \log \log N)$ with exactly $N$ subsequences.
\end{pr}

However, our conjecture still remains unproven. It is not obvious if its potential future proof will be a combinatorial one, or one coming from a seemingly unrelated field -- theory of continued fractions.

\end{document}